\documentclass[preprint,12pt]{elsarticle}

\usepackage{amssymb}

\journal{}

\begin{document}

\begin{frontmatter}



\title{Some Properties of the Coordinate Rings over Finite Fields}


\author{Qinqin Jin,Yongbin Li}

\address{ School of Mathematical Science,University of Electronic Science and Technology of China}

\begin{abstract}
In this paper, we present the  Hilbert's Nullstellensatz in case
of the coordinate rings of a nonempty subset of $S\subseteq K^n$
where $K$ is a finite field $F_q$. Some applications of the
Nullstellensatz  are also discussed.
\end{abstract}

\begin{keyword}
affine variety \sep coordinate ring \sep  Hilbert's
Nullstellensatz

\end{keyword}

\end{frontmatter}


\section{Introduction}
 It is well known that Hilbert's Nullstellensatz plays a
very crucial role in algebraic geometry. It is one of the
milestones in the process of translating algebra into geometry and
geometry into algebra and forms the background for many
applications in algebraic geometry.

Let $ K$ be a  field  and $ K[x_1,\ldots,x_n]$ (or $ K[{\bf x}]$
for short) the ring of polynomials in the variables
$(x_1,\ldots,x_n)$ with coefficients in $ K$.$K^n$ is the
$n$-dimensional affine space over $K$. Let $K\subseteq L$ be a
field extension, $\overline{K}$ be the algebraic closure of $K$.
  For an ideal $ I=<f_{1},\dots,f_{s}>$  generated
by $f_{1},\dots,f_{s}\in K[{\bf x}],$ the set of zeros of $ I$ in
$L^n$ is defined as
$$V_L( I)=\{(a_{1},\dots,a_{n})\in
L^{n}|f_{i}(a_{1},\dots,a_{n})=0,\forall\;1\leq i\leq s\}$$
($V_{\overline{K}}( I)$ is called the affine variety defined by
$I$). It is obvious that $$V_L( I)=V_L( \sqrt{I})$$ where
$\sqrt{I}$ is the radical of $I$, see \cite{Cox, Martin} for the
details. For any set $S\subset K^{n}$, we know that the set
$I(S)=\{f\in K[x_{1},\dots,x_{n}] |  f({\bf a})=0, \forall {\bf
a}\in S\}$ is an ideal in $K[{\bf x}]$.

 {\bf The
 Nullstellensatz}: Let $K$ be an algebraically closed field
and let $I\subset K[{\bf x}]$ be an ideal. Then $I(V_K(I))
=\sqrt{I}.$

In this paper, we will present the following result in case of the
coordinate ring $ K[S]$, which can be identified with $K[{\bf
x}]/I(S)$  where $K=F_{q}$ is a finite field and nonempty set
$S\subseteq K^n$.

{\bf Theorem 1.} Let $K=F_{q}$ be a field where $p$ is a prime
number and $q=p^{k}$ for a positive integer $k$, and a nonempty
subset $S\subset F_{q}^{n}$. If $J$ be a proper ideal of $K[S]$,
then $I_{S}(V_{S}(J))=J$.

 Some applications of the above
theorem are presented in Section 2. Based upon those applications,
several  results presented in \cite{X.G}, which proposes an
efficient characteristic set methods for computing zeros of
polynomial equation systems in a finite field,
 become concise.

\section{Main Result}
Through  this paper, we denote  $S$ by a nonempty subset of
$K^{n}$ where $K$ is a finite field.

{\bf Lemma 2.1.} Let $K$  and $S$ be defined as above, then S is
an affine variety in $K^n$.

{\it Proof:}\/ It is easy to see that $S$ is a finite set since\
$F_{q}^{n}$ is
 finite. Suppose that $S =\{{\bf a}_1,\ldots,{\bf a}_k\}$
 where ${\bf a}_i =(x_{i1},x_{i2},\dots,x_{in})\in F_{q}^{n}$,
 then $S$ is an affine variety since a single point ${\bf a}_i$ is an affine
  variety and the union of finite varieties are also finite variety. $\square$

Based upon the above lemma, with the same notation in \cite{Cox},
we will present the similar
 Nullstellensatz in case of the coordinate ring
$ K[S]\triangleq K[{\bf x}]/I(S)$ where $K=F_{q}$ is a finite
field and nonempty set $S\subseteq K^n$. Note that
$$I(S)=<x_{1}^{q}-x_{1},\ldots,x_{n}^{q}-x_{n}>$$ in case of
$S=F_q^n$, and $K[F_q^n]$ denoted by ${\mathbb R}_q$ in
\cite{X.G}.

 We proceed to discuss some properties of $K[S]$.

{\bf Lemma 2.2.} In the above situation, let $\phi \in K[S]$. Then
$<\phi>=<\phi^{m}>$ for any $m\geq 1$.

{\it Proof:}\/ It is easy to see that$ <\phi^{m}>\subset <\phi> $.
It remains  to show that$<\phi> \subset <\phi^{m}>$. Given any
$\phi \in K[S]$, we have $\phi(\alpha)^{q}-\phi(\alpha)=0$ since
for any $\alpha \in F_{q}^{n}$, $\phi (\alpha)\in F_{q}$, so
$\phi^{q}-\phi=[0]$. Consider that the greatest common divisor of
polynomial $x^{m},x^{q}-x$ is $x$ in $K[x]$, so there exist some
$f(x),g(x)\in K[x]$ such that $f(x)x^{m}+g(x)(x^{q}-x)=x$. It
implies that $$\phi =f(\phi)\phi
^{m}+g(\phi)(\phi^{q}-\phi)=f(\phi)\phi ^{m}.$$ Thus $ <\phi>
\subset <\phi^{m}> $. This completes the proof. $\square$

{\bf Proposition 2.1.} In the above situation, let $J\lhd K[S]$ be
an ideal, then $\sqrt{J}=J$.

 {\it Proof:}\/
For any $\phi\in \sqrt{J}$, there exist $\phi^{m}\in J$ for some
integer $m$, thus $<\phi^{m}>\subset J$. Following from lemma 2.2,
we get $\phi\in J$ and it implies that $ \sqrt{J}\subset J$. The
converse $J\subset \sqrt{J}$ is obvious. This completes the proof.
$\square$

{\bf Remark 2.1:} Proposition 2.1 generalizes Lemma 3 presented in
\cite{X.G}.

In order to explore the correspondence between ideals and
varieties, we introduce the following definitions.

{\bf Definition 2.1.} Let a nonempty set $S\subset K^{n}$ where
$K=F_{q}.$
  For any ideal
 $J=<\phi_{1},\dots ,\phi_{s}> \lhd K[S]$,
 we define
 $$V_{S}(J)=\{(a_{1},\dots ,a_{n})\in S \mid \phi(a_{1},\ldots ,a_{n})=0, \forall \;\phi \:\in J \}.$$
 We call $V_{S}(J)$ a subvariety of $S$. For each nonempty  subset $T\subset S$,
 we define
 $$I_{S}(T)=\{\phi \in K[S]\mid \phi(a_{1},\ldots ,a_{n})=0,\forall \;(a_{1},\ldots ,a_{n})\in T \}.$$

Refer to Exercise 7 in $\S$1 of Chapter 4 in \cite{Cox}, one can
get easily  the following claim. Here, we omit the details.

{\bf Lemma 2.3.} In the above situation, let $J$ be an ideal of $K[S]$.
 If we have $V_{S}(J)=\emptyset$, then there exists a polynomial function $\phi\in J$ such that $V_{S}(<\phi>)=\emptyset$.

{\bf Proposition 2.2~~(Weak Nullstellensatz).} Let $K=F_{q}$ and a
nonempty subset $S\subset F_{q}^{n}$. If $J$ is a proper ideal of
$K[S]$, then $V_{S}(J)\neq \emptyset$.

{\it Proof:}\/
 We will prove this in contrapositive form: if
$V_{S}(J)=\emptyset$, then $J=K[S]$. From lemma $2.3$, we know
that if $V_S(J)=\emptyset$, there exists $\phi \in J$, such that
$V_{S}(<\phi>)=\emptyset$, then $\phi({\bf a})\neq 0$ for any
$\bf{a} \in S$. It implies that $\phi^{q-1}=[1]\in J$, so
$J=K[S]$. This completes the proof. $\square$

Based upon the above results, we proceed to present the following
 Nullstellensatz in the case of the
coordinate ring $K[S].$

{\bf Theorem 2.1~(Hilbert's Nullstellensatz).} Let $K=F_{q}$, a
nonempty subset $S\subset F_{q}^{n}$ and $J$ be
 a proper ideal of $K[S]$. Then $I_{S}(V_{S}(J))=J$.

{\it Proof:}\/ The inclusion $I_{S}(V_{S}(J))\supseteq J$ is
obvious, it remain to show the converse.

 Let
$J=<\phi_{1},\dots,\phi_{s}>\lhd K[S]$, for any $\phi\in
I_{S}(V_{S}(J))$, consider the ideal
$$\widetilde{J}=<\phi_{1},\dots,\phi_{s},[1]-[y]\phi>\subset
K[x_{1},\dots,x_{n},y]/I(\widetilde{S}),$$ where
$$\widetilde{S}=\{(a_{1},\dots,a_{n},a_{n+1})\in F_{q}^{n+1}\mid
(a_{1},\dots,a_{n})\in S,a_{n+1}\in F_{q}\}.$$ We claim that
$V_{\widetilde{S}}(\widetilde{J})=\emptyset$. To see this, let
$(a_{1},\dots,a_{n},a_{n+1})\in \widetilde{S}$.

$(1)$ $(a_{1},\dots,a_{n})$ is a common zero of
 $\phi_{1},\dots ,\phi_{s}$, in this case,
 $$\phi(a_{1},\dots,a_{n})=0$$ since
  $ \phi $ vanishes at any common zero of $\phi_{1},\dots ,\phi_{s}$,
   thus, the polynomial function $[1]-[y]\phi$ takes the value
    $ 1-a_{n+1}\phi(a_{1},\dots,a_{n})\neq 0$
    at the point $(a_{1},\dots,a_{n},a_{n+1})$,
     hence $(a_{1},\dots,a_{n},a_{n+1})$ is not the point of
      $V_{\widetilde{S}}(\widetilde{J})$;

   {$(2)$ $(a_{1},\dots,a_{n})$ is not a common zero of
    $\phi_{1},\dots,\phi_{s}$, in this case, for some $i$ (
     $1\leq i\leq s$), we must have $\phi_{i}(a_{1},\dots,a_{n})\neq 0$.
      Thinking of $\phi_{i}$ as a function of $n+1$ variables
       which does not depend on the last variable, we have
        $\phi_{i}(a_{1},\dots,a_{n+1})\neq 0$,
        so $(a_{1},\dots,a_{n},a_{n+1})$ is not the point of
         $V_{\widetilde{S}}(\widetilde{J})$. This implies
          $V_{\widetilde{S}}(\widetilde{J})=\emptyset$.
          By Proposition 2.2, we conclude that $[1]\in \widetilde{J}$.
           Thus, $$[1]=\sum_{i=1}^np_{i}\phi_{i}+q([1]-[y]\phi)$$ for
            some functions  $p_{i},q \in K[x_{1},\dots,x_{n},y]/I(\widetilde{S})$.
           Now set $[y]=[1]/\phi$, and multiply both side by a power of $\phi^{m}$,
           where $m$ is large enough to clear all the denominators.
           This yields $\phi^{m}=\sum_{i=1}^{s}h_{i}\phi_{i}$. So $\phi \in \sqrt{J}$ and we get $I_S(V_{S}(J))\subseteq \sqrt{J}$.}

Since $\sqrt{J}=J$ by Lemma 2.1, we have that $I_{S}(V_{S}(J))=J$.
This completes the proof. $\square$

 {\bf  Corollary  2.3.} Let
$K=F_{q},S=F_{q}^{n}$ and $[f]\in K[S]$,
$$V_{S}(<[f]>)=F_{q}^{n}$$ if and only if $[f]\equiv [0]$.

{\it Proof:}\/
 If $[f]\equiv [0]$, then obviously $V_{S}(<[f]>)=F_{q}^{n}$.
Conversely, since $$I_{S}(F_{q}^{n})=\{[g]\in K[F_{q}^{n}] \mid
g({\bf a})=0,\forall {\bf a} \in F_{q}^{n}\}=<[0]>,$$ this shows
$I_{S}(V_{S}(<[f]>))=<[0]>$. It follows from Theorem 2.1 that
$[f]\equiv [0].$ This completes the proof. $\square$

{\bf Remark 2.2:} Lemmas 4 and 5  presented in \cite{X.G} are
direct consequences of a special case of Theorem 2.1 when
$S=F_q^n.$

\section{Application}

 As a consequence of Theorem 2.1 and Remark 2.2,
  the proof of the  following result is similar to the
   argument presented in
  Theorem 5
of $\S 4$ of Chapter 5 in \cite{Cox}, here we omit the details.

{\bf Proposition 3.1.} Let $K=F_{q}$, nonempty subset $S\subset
F_{q}^{n}$.

(1) The correspondences
$$I_S:P(S)\rightarrow \{J\lhd K[S]\},\;V_S:\{J\lhd K[S]\}\rightarrow P(S) $$
where $P(S)$ is the power set of $S$ are inclusion-reversing
bijections and are inverses of each other.

(2) Under the correspondence given in (1), points of $S$
correspond to maximal ideals of $K[S]$.

 The following consequences are helpful to better understand the
coordinate ring $ K[S]$ and some results presented in \cite{X.G}.

{\bf Proposition 3.2.} Let $I,J$ be ideals in $K[S]$, then we have
$V_{S}(I:J)=V_{S}(I) \backslash V_{S}(J)$.

 {\it Proof:}\/
 We claim that $V_{S}(I:J) \supset V_{S}(I) \setminus
V_{S}(J)$. The proof is similar to the argument given in Theorem 7
of $\S 4$ of Chapter 4 in \cite{Cox}, we omit the details. Now we
prove the opposite inclusion. Let ${\bf a} \in V_{S}(I:J)$, that
is to say, if $\phi \psi \in I$ for all $\psi \in J$, then
$\phi({\bf a})=0$. Let $\phi \in I_{S}(V_{S}(I) \setminus
V_{S}(J))$, if $\psi \in J$, then $\phi \psi$ vanishes on
$V_{S}(I)\setminus V_{S}(J)$ and $\psi$ on $V_{S}(J)$. Thus, by
the theorem $2.1$, $\phi \psi \in I$ and $\phi({\bf a})=0$. Hence,
${\bf a} \in V_{S}(I_{S}(V_{S}(I)\setminus V_{S}(J)))$. This
establishes that $V_{S}(I:J)\subset V_{S}(I_{S}(V_{S}(I)\setminus
V_{S}(J)))$, and completes the proof. $\square$

{\bf Corollary 3.1.} Let $T_{1},T_{2} \in P(S) $, then $T_{1}
\backslash T_{2}=V_{S}(I_{S}(T_{1}) \backslash I_{S}(T_{2}))$, and
$$I_{S}(S_{1}):I_{S}(S_{2})=I_{S}(S_{1}\setminus S_{2}).$$

The next result  extends Lemma 7 in \cite{X.G} to the case of when
any nonempty set $S\subseteq F_{q}^{n}$.

 {\bf Proposition 3.3.} In the above situation, let $[f],[g]$ and $[h]$ be
polynomial
functions in $K[S]$. We have\\
(1) $V_{S}(<[f]^{q-1}[g]^{q-1}-[1]>)=V_{S}(<[f]^{q-1}-[1]>)\cap
V_{S}(<[g]^{q-1}-[1]>)$, and
$$<f^{q-1}g^{q-1}-1>=<f^{q-1}-1>+<g^{q-1}-1>;$$
(2) $V_{S}(<[f]^{q-1}[g]^{q-1}-[f]^{q-1}-[g]^{q-1}>)=V_{S}(<[f]>)\cap V_{S}(<[g]>)$, and
$$<[f]^{q-1}[g]^{q-1}-[f]^{q-1}-[g]^{q-1}>=<[f]>+<[g]>.$$

{\it Proof:}\/ (1) Let ${\bf a} \in V_{S}(<[f]^{q-1}-[1])\cap
V_{S}([g]^{q-1}-[1]>)$,
 then $f({\bf a})^{q-1}=1$ and $g({\bf a})^{q-1}=1$,
 so we have $[f]^{q-1}[g]^{q-1}-[1]\equiv [0]$ and
 $${\bf a} \in V_{S}(<[f]^{q-1}[g]^{q-1}-[1]>).$$
Conversely, suppose ${\bf a} \in V_{S}(<[f]^{q-1}[g]^{q-1}-[1]>)$.
Since $f({\bf a}),g({\bf a})\in F_{q}$, $f({\bf a})=0$ or $f({\bf
a})^{q-1}=1$, the same with $g({\bf a})$. Only when
 $f({\bf a})^{q-1}=1$ and $g({\bf a})^{q-1}=1$, we have
 $[f]^{q-1}[g]^{q-1}=[1],$ so
 $${\bf a} \in V_{S}(<[f]^{q-1}-[1])\cap V_{S}([g]^{q-1}-[1])>.$$ It
 follows from Theorem 2.1 that
 \begin{eqnarray*}
<f^{q-1}g^{q-1}-1>&=&I_S(V_S(<f^{q-1}g^{q-1}-1))\\
 &=&I_S(V_{S}(<[f]^{q-1}-[1]>)\cap V_{S}(<[g]^{q-1}-[1]>))\\
 &=&I_S(V_{S}(<[f]^{q-1}-[1]>))+I_S( V_{S}(<[g]^{q-1}-[1]>))\\
 &=&<[f]^{q-1}-[1]>+<[g]^{q-1}-[1]>.
\end{eqnarray*}

(2) The proof of part (2) is similar to the proof part (1) and
will be omitted. This completes the proof. $\square$





\label{}
\bibliographystyle{elsarticle-num}
\bibliography{<your-bib-database>}


\end{document}